\documentclass[12pt]{article}
\usepackage{mathptmx}
\usepackage{graphicx}
\usepackage{bm}
\usepackage{amssymb}
\makeatletter
\author{Robert P. C. de Marrais \footnote{Email address:  rdemarrais@alum.mit.edu} \\ \noindent
\emph{Thothic Technology Partners, P.O.Box 3083, Plymouth MA 02361}}

\title{Placeholder Substructures III:~~  A Bit-String-Driven ``Recipe Theory'' for Infinite-Dimensional Zero-Divisor Spaces}

\begin{document}
\maketitle
\makeatother

\begin{abstract}
Zero-divisors (ZDs) derived by Cayley-Dickson Process (CDP) from
N-dimensional hypercomplex numbers ($N$ a power of $2$, and at least
$4$) can represent singularities and, as $N \rightarrow \infty$,
fractals -- and thereby, scale-free networks.  Any integer $> 8$ and
not a power of $2$ generates a meta-fractal or \textit{Sky} when it
is interpreted as the \textit{strut constant} (S) of an ensemble of
octahedral vertex figures called \textit{Box-Kites} (the fundamental
ZD building blocks). Remarkably simple bit-manipulation rules or
\textit{recipes} provide tools for transforming one fractal genus
into others within the context of Wolfram's Class 4 complexity.
\end{abstract}

\section{The Argument So Far}
In Parts I[1] and II[2], the basic facts concerning zero-divisors
(ZDs) as they arise in the hypercomplex context were presented and
proved. ``Basic,'' in the context of this monograph, means seven
things. First, they emerged as a side-effect of applying CDP a
minimum of 4 times to the Real Number Line, doubling dimension to
the Complex Plane, Quaternion 4-Space, Octonion 8-Space, and 16-D
Sedenions. With each such doubling, new properties were found:  as
the price of sacrificing counting order, the Imaginaries made a
general theory of equations and solution-spaces possible; the
non-commutative nature of Quaternions mapped onto the realities of
the manner in which forces deploy in the real world, and led to
vector calculus; the non-associative nature of Octonions, meanwhile,
has only come into its own with the need for necessarily
unobservable quantities (because of conformal field-theoretical
constraints)in String Theory.  In  the Sedenions, however, the most
basic assumptions of all -- well-defined notions of field and
algebraic norm (and, therefore, measurement) -- break down, as the
phenomena correlated with their absence, zero-divisors, appear
onstage (never to leave it for all higher CDP dimension-doublings).

Second thing:  ZDs require at least two differently-indexed
imaginary units to be defined, the index being an integer larger
than 0 (the CDP index of the Real Unit) and less than $2^{N}$ for a
given CDP-generated collection of $2^{N}$-ions.  In ``pure CDP,''
the enormous number of alternative labeling schemes possible in any
given $2^{N}$-ion level are drastically reduced by assuming that
units with such indices interact by XOR-ing:  the index of the
product of any two is the XOR of their indices.  Signing is more
tricky; but, when CDP is reduced to a 2-rule construction kit, it
becomes easy:  for index $u < {\bf G}$, ${\bf G}$ the Generator of
the $2^{N}$-ions (i.e., the power of 2 immediately larger than the
highest index of the predecessor $2^{N-1}$-ions), Rule 1 says $i_{u}
\cdot i_{\bf G} = + i_{(u + {\bf G})}$.   Rule 2 says take an
associative triplet $(a, b, c)$, assumed written in CPO (short for
``cyclically positive order'':  to wit, $a \cdot b = +c$, $b \cdot c
= +a$, and $c \cdot a = + b$).  Consider, for instance, any $(u,
{\bf G}, {\bf G} + u)$ index set.  Then three more such associative
triplets (henceforth, \textit{trips}) can be generated by adding
{\bf G} to two of the three, then switching their resultants' places
in the CPO scheme. Hence, starting with the Quaternions' $(1, 2, 3)$
(which we'll call a \textit{Rule 0} trip, as it's inherited from a
prior level of CDP induction), Rule 1 gives us the trips $(1, 4,
5)$, $(2, 4, 6)$, and $(3, 4, 7)$, while Rule 2 yields up the other
4 trips defining the Octonions:  $(1, 7, 6)$, $(2, 5, 7)$, and $(3,
6, 5)$.  Any ZD in a given level of $2^{N}$-ions will then have
units with one index $< {\bf G}$, written in lowercase, and the
other index $> {\bf G}$, written in uppercase. Such pairs,
alternately called ``dyads'' or ``Assessors,'' saturate the diagonal
lines of their planes, which diagonals never mutually zero-divide
each other (or make \textit{DMZs}, for "divisors (or dyads) making
zero"), but only make DMZs with other such diagonals, in other such
Assessors. (This is, of course, the opposite situation from the
projection operators of quantum mechanics, which are diagonals in
the planes formed by Reals and dimensions spanned by Pauli spin
operators contained within the 4-space created by the Cartesian
product of two standard imaginaries.)

Third thing:  Such ZDs are not the only possible in CDP spaces; but
they define the ``primitive'' variety from which ZD spaces
saturating more than 1-D regions can be articulated.  A not quite
complete catalog of these can be found in our first monograph on the
theme [3]; a critical kind which was overlooked there, involving the
Reals (and hence, providing the backdrop from which to see the
projection-operator kind as a degenerate type), were first discussed
more recently [4]. (Ironically, these latter are the easiest sorts
of composites to derive of any:  place the two diagonals of a DMZ
pairing with differing internal signing on axes of the same plane,
and consider the diagonals \textit{they} make with each other!)  All
the primitive ZDs in the Sedenions can be collected on the vertices
of one of 7 copies of an Octahedron in the \textit{Box-Kite}
representation, each of whose 12 edges indicates a two-way ``DMZ
pathway,'' evenly divided between 2 varieties.  For any vertex V,
and $k$ any real scalar, indicate the diagonals this way:
$\verb|(V,/)|  = k \cdot (i_{v} + i_{V})$, while $\verb|(V, \)| = k
\cdot (i_{v} - i_{V})$. 6 edges on  a Box-Kite will always have
\textit{negative edge-sign} (with \textit{unmarked} ET cell entries:
 see the ``sixth thing''). For vertices M and N, exactly two
DMZs run along the edge joining them, written thus:

\begin{center}
$\verb|(M,/)| \cdot \verb|(N, \)|$ $= \verb|(M, \)|$ $\cdot$
$\verb|(N,/)|$ $ = 0$
\end{center}

The other 6 all have \textit{positive edge-sign}, the diagonals of
their two DMZs having same slope (and \textit{marked} -- with
leading dashes -- ET cell entries):

\begin{center}
$\verb|(Z,/)| \cdot \verb|(V, /)| = $ $\verb|(Z, \)| \cdot
\verb|(V,\)| = 0$
\end{center}

Fourth thing:  The edges always cluster similarly, with two opposite
faces among the 8 triangles on the Box-Kite being spanned by 3
negative edges (conventionally painted red in color renderings),
with all other edges being positive (painted blue).  One of the red
triangles has its vertices' 3 low-index units forming a trip;
writing their vertex labels conventionally as A, B, C, we find there
are in fact always 4 such trips cycling among them:  $(a, b, c)$,
the \textit{L-trip}; and the three \textit{U-trips} obtained by
replacing all but one of the lowercase labels in the L-trip with
uppercase:  $(a, B, C)$; $(A, b, C)$; $(A, B, c)$.  Such a 4-trip
structure is called a \textit{Sail}, and a Box-Kite has 4 of them:
the \textit{Zigzag}, with all negative edges, and the 3
\textit{Trefoils}, each containing two positive edges extending from
one of the Zigzag vertices to the two vertices opposite its
\textit{Sailing partners}.  These opposite vertices are always
joined by one of the 3 negative edges comprising the Vent which is
the Zigzag's opposite face.  Again by convention, the vertices
opposite A, B, C are written F, E, D in that order; hence, the
Trefoil Sails are written $(A, D, E)$; $(F, D, B)$, and $(F, C, E)$,
ordered so that their lowercase renderings are equivalent to their
CPO L-trips.  The graphical convention is to show the Sails as
filled in, while the other 4 faces, like the Vent, are left empty:
they show ``where the wind blows'' that keeps the Box-Kite aloft. A
real-world Box-Kite, meanwhile, would be held together by 3 dowels
(of wood or plastic, say) spanning the joins between the only
vertices left unconnected in our Octahedral rendering: the
\textit{Struts} linking the \textit{strut-opposite} vertices (A, F);
(B, E); (C, D).

Fifth thing:  In the Sedenions, the 7 isomorphic Box-Kites are
differentiated by which Octonion index is missing from the vertices,
and this index is designated by the letter {\bf S}, for
``signature,'' ``suppressed index,'' or \textit{strut constant}.
This last designation derives from the invariant relationship
obtaining in a given Box-Kite between {\bf S} and the indices in the
Vent and Zigzag termini (V and Z respectively) of any of the 3
Struts, which we call the ``First Vizier'' or VZ1.  This is one of 3
rules, involving the three Sedenion indices always missing from a
Box-Kite's vertices: ${\bf G}$, ${\bf S}$, and their simple sum
${\bf X}$ (which is also their XOR product, since ${\bf G}$ is
always to the left of the left-most bit in ${\bf S}$).  The Second
Vizier tells us that the L-index of either terminus with the U-index
of the other always form a trip with ${\bf G}$, and it true as
written for all $2^{N}$-ions.  The Third shows the relationship
between the L- and U- indices of a given Assessor, which always form
a trip with ${\bf X}$.  Like the First, it is true as written only
in the Sedenions, but as an unsigned statement about indices only,
it is true universally.  (For that reason, references to VZ1 and VZ3
hereinout will be assumed to refer to the \textit{unsigned}
versions.)  First derived in the last section of Part I, reprised in
the intro of Part II, we write them out now for the third and final
time in this monograph:

\begin{center}
VZ1:  $v \cdot z = V \cdot Z = {\bf S}$

\smallskip
VZ2:  $Z \cdot v = V \cdot z = {\bf G}$

\smallskip
VZ3:  $V \cdot v = z \cdot Z = {\bf X}$.
\end{center}

Rules 1 and 2, the Three Viziers, plus the standard Octonion
labeling scheme derived from the simplest finite projective group,
usually written as PSL(2,7), provide the basis of our toolkit. This
last becomes powerful due to its capacity for recursive re-use at
all levels of CDP generation, not just the Octonions.  The simplest
way to see this comes from placing the unique Rule 0 trip provided
by the Quaternions on the circle joining the 3 sides' midpoints,
with the Octonion Generator's index, 4, being placed in the center.
Then the 3 lines leading from the Rule 0 trip's (1, 2, 3) midpoints
to their opposite angles -- placed conventionally in clockwise order
in the midpoints of the left, right, and bottom sides of a triangle
whose apex is at 12 o'clock -- are CPO trips forming the Struts,
while the 3 sides themselves are the Rule 2 trips.  These 3 form the
L-index sets of the Trefoil Sails, while the Rule 0 trip provides
the same service for the Zigzag.  By a process analogized to tugging
on a slipcover (Part I) and pushing things into the central zone of
hot oil while wok-cooking (Part II), all 7 possible values of {\bf
S} in the Sedenions, not just the 4, can be moved into the center
while keeping orientations along all 7 lines of the Triangle
unchanged. Part II's critical Roundabout Theorem tells us, moreover,
that all $2^{N}$-ion ZDs, for all $N
> 3$, are contained in Box-Kites as their minimal ensemble size.
Hence, by placing the appropriate ${\bf G}$, ${\bf S}$, or ${\bf X}$
in the center of a PSL(2,7) triangle, with a suitable Rule 0 trip's
indices populating the circle, any and all \textit{candidate}
primitive ZDs can be discovered and situated.

Sixth thing:  The word ``candidate'' in the above is critical; its
exploration was the focus of Part II. For, starting with $N = 5$ and
hence ${\bf G = 16}$ (which is to say, in the 32-D Pathions), whole
Box-Kites can be suppressed (meaning, all 12 edges, and not just the
Struts, no longer serve as DMZ pathways). But for all $N$, the full
set of candidate Box-Kites are viable when ${\bf S} \leq 8$ or equal
to some higher power of 2. For all other ${\bf S}$ values, though,
the phenomenon of \textit{carrybit overflow} intervenes -- leading,
ultimately, to the ``meta-fractal'' behavior claimed in our
abstract.  To see this, we need another mode of representation, less
tied to 3-D visualizing, than the Box-Kite can provide.  The answer
is a matrix-like method of tabulating the products of candidate ZDs
with each other, called \textit{Emanation Tables} or \textit{ETs}.
The L-indices only of all candidate ZDs are all we need indicate
(the U-indices being forced once ${\bf G}$ is specified); these will
saturate the list of allowed indices $< {\bf G}$, save for the value
of ${\bf S}$ whose choice, along with that of ${\bf G}$, fixes an
ET.  Hence, the unique ET for given ${\bf G}$ and ${\bf S}$ will
fill a square spreadsheet whose edge has length $2^{N-1}- 2$.
Moreover, a cell entry (r,c) is only filled when row and column
labels R and C form a DMZ, which can never be the case along an ET's
long diagonals:  for the diagonal starting in the upper left corner,
R xor R = 0, and the two diagonals within the same Assessor, can
never zero-divide each other; for the righthand diagonal, the
convention for ordering the labels (ascending counting order from
the left and top, with any such label's strut-opposite index
immediately being entered in the mirror-opposite positions on the
right and bottom) makes R and C strut-opposites, hence also unable
to form DMZs.

For the Sedenions, we get a 6 x 6 table, 12 of whose cells (those on
long diagonals) are empty:  the 24 filled cells, then, correspond to
the two-way traffic of ``edge-currents'' one imagines flowing
between vertices on a Box-Kite's 12 edges.  A computational
corollary to the Roundabout Theorem, dubbed the \textit{Trip-Count
Two-Step}, is of seminal importance.  It connects this most basic
theorem of ETs to the most basic fact of associative triplets,
indicated in the opening pages of Part I, namely: for any N, the
number $Trip_{N}$ of associative triplets is found, by simple
combinatorics, to be $(2^{N} - 1)(2^{N} - 2)/3!$ -- 35 for the
Sedenions, 155 for the Pathions, and so on. But, by Trip-Count
Two-Step, we also know that \textit{the maximum number of Box-Kites
that can fill a $2^{N}$-ion ET = $Trip_{N-2}$.} For ${\bf S}$ a
power of 2, beginning in the Pathions (for ${\bf S} = 2^{5 - 2} =
8$), the Number Hub Theorem says the upper left quadrant of the ET
is an unsigned multiplication table of the $2^{N-2}$-ions in
question, with the 0's of the long diagonal (indicated Real negative
units) replaced by blanks -- a result effectively synonymous with
the Trip-Count Two-Step.

Seventh thing:  We found, as Part II's argument wound down, that the
2 classes of ETs found in the Pathions -- the ``normal'' for ${\bf
S} \leq 8$, filled with indices for all 7 possible Box-Kites, and
the ``sparse'' so-called Sand Mandalas, showing only 3 Box-Kites
when $8 < {\bf S} < 16$, were just the beginning of the story.  A
simple formula involving just the bit-string of ${\bf s}$ and ${\bf
g}$, where the lowercase indicates the values of ${\bf S}$ and ${\bf
G}$ modulo ${\bf G}/2$, gave the prototype of our first
\textit{recipe}: all and only cells with labels R or C, or content P
( = R xor C ), are filled in the ET.  The 4 ``missing Box-Kites''
were those whose L-index trip would have been that of a Sail in the
$2^{N-1}$ realm with ${\bf S} = {\bf s}$ and ${\bf G} = {\bf g}$.
The sequence of 7 ETs, viewed in ${\bf S}$-increasing succession,
had an obvious visual logic leading to their being dubbed a
\textit{flip-book}. These 7 were obviously indistinguishable from
many vantages, hence formed a \textit{spectrographic band}.  There
were 3 distinct such bands, though, each typified by a Box-Kite
count common to all band-members, demonstrable in the ETs for the
64-D Chingons.  Each band contained ${\bf S}$ values bracketed by
multiples of 8 (either less than or equal to the higher, depending
upon whether the latter was or wasn't a power of 2).  These were
claimed to underwrite behaviors in all higher $2^{N}$-ion ETs,
according to 3 rough patterns in need of algorithmic refining in
this Part III. Corresponding to the first unfilled band, with ETs
always missing $4^{N-4}$ of their candidate Box-Kites for $N > 4$,
we spoke of \textit{recursivity}, meaning the ETs for constant ${\bf
S}$ and increasing $N$ would all obey the same recipe, properly
abstracted from that just cited above, empirically found among the
Pathions for ${\bf S} > 8$. The second and third behaviors, dubbed,
for ${\bf S}$ ascending, \textit{(s,g)-modularity} and
\textit{hide/fill involution} respectively, make their first
showings in the Chingons, in the bands where $16 < {\bf S} \leq 24$,
and then where $24 < {\bf S} < 32$.  In all such cases, we are
concerned with seeing the ``period-doubling'' inherent in CDP and
Chaotic attractors both become manifest in a repeated doubling of ET
edge-size, leading to the fixed-${\bf S}$, $N$ increasing analog of
the fixed-$N, {\bf S}$ increasing flip-books first observed in the
Pathions, which we call \textit{balloon-rides}.  Specifying and
proving their workings, and combining all 3 of the above-designated
behaviors into the ``fundamental theorem of zero-division algebra,''
will be our goals in this final Part III.  Anyone who has read this
far is encouraged to bring up the graphical complement to this
monograph, the 78-slide Powerpoint show presented at NKS 2006 [5],
in another window. (Slides will be referenced by number in what
follows.)

\section{$8 <{\bf S} < 16, N \rightarrow \infty$ :  Recursive Balloon Rides in the Whorfian Sky}
We know that any ET for the $2^{N}$-ions is a square whose edge is
$2^{N-1} - 2$ cells.  How, then, can any simply recursive rule
govern exporting the structure of one such box to analogous boxes
for progressively higher $N$?  The answer: \textit{include the label
lines} -- not just the column and row headers running across the top
and left margins, but their strut-opposite values, placed along the
bottom and right margins, which are mirror-reversed copies of the
label-lines \textit{(LLs)} proper to which they are parallel. This
increases the edge-size of the ET box to $2^{N-1}$.

\medskip

\noindent {\small Theorem 11.}  For any fixed ${\bf S} > 8$ and not
a power of $2$, the row and column indices comprising the Label
Lines (LLs) run along the left and top borders of the $2^{N}$-ion ET
"spreadsheet" for that ${\bf S}$.  Treat them as included in the
spreadsheet, \textit{as labels}, by adding a row and column to the
given square of cells, of edge $2^{N-1} - 2$, which comprises the ET
proper. Then add another row and column to include the
strut-opposite values of these labels' indices in ``mirror LLs,''
running along the opposite edges of a now $2^{N-1}$-edge-length box,
whose four corner cells, like the long diagonals they extend, are
empty. When, for such a fixed ${\bf S}$, the ET for the
$2^{N+1}$-ions is produced, the values of the 4 sets of LL indices,
bounding the contained $2^{N}$-ion ET, correspond, \textit{as cell
values}, to actual DMZ P-values in the bigger ET, residing in the
rows and columns labeled by the contained ET's ${\bf G}$ and ${\bf
X}$ (the containing ET's $g$ and $g + {\bf S}$).  Moreover, all
cells contained in the box they bound in the containing ET have
P-values (else blanks) exactly corresponding to -- and
\textit{including edge-sign markings} of -- the positionally
identical cells in the $2^{N}$-ion ET:  those, that is, for which
the LLs act as labels.

\medskip

\noindent \textit{Proof.}  For all strut constants of interest,
${\bf S} < g ( = {\bf G/2})$; hence, all labels up to and including
that immediately adjoining its own strut constant (that is, the
first half of them) will have indices monotonically increasing, up
to and at least including the midline bound, from $1$ to $g - 1$.
When $N$ is incremented by $1$, the row and column midlines
separating adjoining strut-opposites will be cut and pulled apart,
making room for the labels for the $2^{N+1}$-ion ET for same ${\bf
S}$, which middle range of label indices will also monotonically
increase, this time from the current $2^N$-ion generation's $g$ (and
prior generation's ${\bf G}$), up to and at least including its own
midline bound, which will be $g$ plus the number of cells in the LL
inherited from the prior generation, or $g/2 - 1$. The LLs are
therefore contained in the rows and columns headed by $g$ and its
strut opposite, $g + {\bf S}$.  To say that the immediately prior
CDP generation's ET labels are converted to the current generation's
P-values in the just-specified rows and columns is equivalent to
asserting the truth of the following calculation:

\smallskip

\begin{center}
$(g + u) + (sg) \cdot ({\bf G} + g + u_{opp})$

\underline{$\;\;\;\;\;\;g   \;\;\;\;\; + \;\;\;\;\;  ({\bf G} + g +
{\bf S})\;\;\;\;\;\;$}
\smallskip

$- (vz)\cdot ({\bf G} + u_{opp}) \;\;\; + (vz) \cdot (sg) \cdot u$

\underline{$ + u  \;\;\;\;\; - \;\;\;\;\; (sg) \cdot ({\bf G} +
u_{opp})$}

\smallskip
$0$ only if $vz = (-sg)$

\end{center}

\smallskip

Here, we use two binary variables, the inner-sign-setting $sg$, and
the Vent-or-Zigzag test, based on the First Vizier.  Using the two
in tandem lets us handle the normal and ``Type II'' box-kites in the
same proof.  Recall (and see Appendix B of Part II for a quick
refresher) that while the ``Type I'' is the only type we find in the
Sedenions, we find that a second variety emerges in the Pathions,
indistinguishable from Type I in most contexts of interest to us
here:  the orientation of 2 of the 3 struts will be reversed (which
is why VZ1 and VZ3 are only true generally when unsigned). For a
Type I, since ${\bf S} < g$, we know by Rule 1 that we have the trip
$({\bf S}, g, g + {\bf S})$; hence, $g$ -- for all $2^{N}$-ions
beyond the Pathions, where the Sand Mandalas' $g = 8$ is the L-index
of the Zigzag B Assessor -- must be a Vent (and its strut-opposite,
$g + {\bf S}$, a Zigzag). For a Type II, however, this is
necessarily so only for 1 of the 3 struts -- which means, per the
equation above, that sg must be reversed to obtain the same result.
Said another way, we are free to assume either signing of $vz$ means
+1, so the ``only if'' qualifying the zero result is informative.
 It is $u$ and its relationship to $g + u$ that is of interest
here, and this formulation makes it easier to see that the products
hold for arbitrary LL indices $u$ \textit{or} their strut-opposites.
But for this, the term-by-term computations should seem routine: the
left bottom is the Rule 1 outcome of $(u, g, g+u)$:  obviously, any
$u$ index must be less than $g$.  To its right, we use the trip
$(u_{opp}, g, g + u_{opp}) \rightarrow ({\bf G} + g + u_{opp}, g,
{\bf G} + u_{opp})$, whose CPO order is opposite that of the
multiplication.  For the top left, we use $(u, {\bf S}, u_{opp})$ as
limned above, then augment by $g$, then ${\bf G}$, leaving $u_{opp}$
unaffected in the first augmenting, and $g + u$ in the second.
Finally, the top right (ignoring $sg$ and $vz$ momentarily) is
obtained this way: $(u,{\bf S}, u_{opp}) \rightarrow (u, g +
u_{opp}, g + {\bf S}) \rightarrow (u, {\bf G} + g + s, {\bf G} + g +
u_{opp})$; ergo, $+u$.

Note that we cannot eke out any information about edge-sign marks
from this setup:  since labels, as such, have no marks, we have
nothing to go on -- unlike all other cells which our recursive
operations will work on.  Indeed, the exact algorithmic
determination of edge-sign marks for labels is not so trivial:  as
one iterates through higher $N$ values, some segments of LL indexing
will display reversals of marks found in the ascending or descending
left midline column, while other segments will show them unchanged
-- with key values at the beginnings and ends of such octaves
(multiples of $8$, and sums of such multiples with ${\bf S}\; mod \;
8$) sometimes being reversed or kept the same irrespective of the
behavior of the terms they bound.  Fortunately, such behaviors are
of no real concern here -- but they are, nevertheless, worth
pointing out, given the easy predictability of other edge-sign marks
in our recursion operations.

Now for the ET box within the labels:  if all values (including
edge-sign marks) remain unchanged as we move from the $2^{N}$-ion ET
to that for the $2^{N+1}$-ions, then one of 3 situations must
obtain:  the inner-box cells have labels $u, v$ which belong to some
Zigzag L-trip $(u, v, w)$; or, on the contrary, they correspond to
Vent L-indices -- the first two terms in the CPO triplet $(w_{opp},
v_{opp}, u)$, for instance; else, finally, one term is a Vent, the
other a Zigzag (so that inner-signs of their multiplied dyads are
both positive):  we will write them, in CPO order, $v_{opp}$ and
$u$, with third trip member $w_{opp}$.  Clearly, we want all the
products in the containing ET to indicate DMZs only if the inner
ET's cells do similarly.  This is easily arranged:  for the
containing ET's cells have indices identical to those of the
contained ET's, save for the appending of $g$ to both (and ditto for
the U-indices).

\medskip

\textit{Case 1:}  If $(u, v, w)$ form a Zigzag L-index set, then so
do $(g + v, g + u, w)$, so markings remain unchanged; and if the
$(u,v)$ cell entry is blank in the contained, so will be that for
$(g + u, g + v)$ in its container. In other words, the following
holds:

\smallskip

\begin{center}
$(g + v) + (sg) \cdot ({\bf G} + g + v_{opp})$

\underline{$\;\;\;(g + u)   \;\;\;\; + \;\;\;({\bf G} + g +
u_{opp})\;\;\;$}
\smallskip

$- ({\bf G} + w_{opp}) \;\;\; - (sg) \cdot w$

\underline{$ - w \;\;\;\;\; - (sg) \cdot ({\bf G} + w_{opp})$}

\smallskip
$0$ only if $sg = (-1)$

\end{center}

\pagebreak

$(g + u) \cdot (g + v) = P:$  $(u,v,w) \rightarrow (g + v, g + u,
w)$; hence, $(- w)$.

$(g + u) \cdot (sg) \cdot ({\bf G} + g + v_{opp}) = P:$  $(u,
w_{opp}, v_{opp})$ $\rightarrow (g + v_{opp}, w_{opp}, g + u)$
$\rightarrow ({\bf G} + w_{opp}, {\bf G} + g + v_{opp}, g + u)$;
hence, $(sg) \cdot ( - ({\bf G} + w_{opp}))$.

$({\bf G} + g + u_{opp}) \cdot (g + v) = P:$  $(u_{opp}, w_{opp},
v)$ $\rightarrow (g + v, w_{opp}, g + u_{opp})$ $\rightarrow (g + v,
{\bf G} + g + u_{opp}, {\bf G} + w_{opp})$; hence, $(- ({\bf G} +
w_{opp}))$.

$({\bf G} + g + u_{opp}) \cdot ({\bf G} + g + v_{opp}) = P:$ Rule 2
twice to the same two terms yields the same result as the terms in
the raw, hence $(- w)$.

Clearly, cycling through $(u,v,w)$ to consider $(g + v) \cdot (g +
w)$ will give the exactly analogous result, forcing two (hence
three) negative inner-signs in the candidate Sail; hence, if we have
DMZs at all, we have a Zigzag Sail.

\medskip

\textit{Case 2:}  The product of two Vents must have negative
edge-sign, and there's no cycling through same-inner-signed products
as with the Zigzag, so we'll just write our setup as a one-off, with
upper inner-sign explicitly negative, and claim its outcome true.

\smallskip

\begin{center}
$(g + v_{opp}) - ({\bf G} + g + v)$

\underline{$\;(g + w_{opp})   \; + \;({\bf G} + g + w)\;$}
\smallskip

$+ ({\bf G} + u_{opp}) \;\;\;\;\; + u$

\underline{$ - u  \;\;\;\;\; - ({\bf G} + u_{opp})$}

\smallskip
$0$

\end{center}

\smallskip

$(g + w_{opp}) \cdot (g + v_{opp}) = P:$  $(w_{opp}, v_{opp}, u)$
$\rightarrow (g + v_{opp}, g + w_{opp}, u)$; hence, $(- u)$.

$(g + w_{opp}) \cdot ({\bf G} + g + v) = P:$  $(w_{opp}, v,
u_{opp})$ $\rightarrow (g + v, g + w_{opp}, u_{opp})$ $\rightarrow
({\bf G} + u_{opp}, g + w_{opp}, {\bf G} + g + v)$; but inner sign
of upper dyad is negative, so $(- ({\bf G} + u_{opp}))$.

$({\bf G} + g + w) \cdot (g + v_{opp}) = P:$  $(v_{opp}, u_{opp},
w)$ $\rightarrow (g + w, u_{opp}, g + v_{opp})$ $\rightarrow ({\bf
G} + u_{opp}, {\bf G} + g + w, g + v_{opp})$; hence, $(+ ({\bf G} +
u_{opp}))$.

$({\bf G} + g + w) \cdot ({\bf G} + g + v) = P:$ Rule 2 twice to the
same two terms yields the same result as the terms in the raw; but
inner sign of upper dyad is negative, so $(+ u)$.

\medskip

\textit{Case 3:}  The product of Vent and Zigzag displays same inner
sign in both dyads; hence the following arithmetic holds:

\pagebreak

\begin{center}
$(g + u) + ({\bf G} + g + u_{opp})$

\underline{$\; (g + v_{opp})   \; + \;({\bf G} + g + v)\;$}
\smallskip

$- ({\bf G} + w) \;\;\;\;\; + w_{opp}$

\underline{$ - w_{opp}  \;\;\;\;\; + ({\bf G} + w)$}

\smallskip
$0$

\end{center}

\smallskip

The calculations are sufficiently similar to the two prior cases as
to make their writing out tedious.  It is clear that, in each of our
three cases, content and marking of each cell in the contained ET
and the overlapping portion of the container ET are identical. $\;\;
\blacksquare$

\medskip

To highlight the rather magical label/content involution that occurs
when $N$ is in- or de- cremented, graphical realizations of such
nested patterns, as in Slides 60-61, paint LLs (and labels proper) a
sky-blue color.  The bottom-most ET being overlaid in the central
box has $g =$ the maximum high-bit in ${\bf S}$, and is dubbed the
\textit{inner skybox}.  The degree of nesting is strictly measured
by counting the number of bits $B$ that a given skybox's $g$ is to
the left of this strut-constant high-bit.  If we partition the inner
skybox into quadrants defined by the midlines, and count the number
$Q$ of quadrant-sized boxes along one or the other long diagonal, it
is obvious that the inner skybox itself has $B = 0$ and $Q = 1$; the
nested skyboxes containing it have $Q = 2^{B}$.  If recursion of
skybox nesting be continued indefinitely -- to the fractal limit,
which terminology we will clarify shortly -- the indices contained
in filled cells of any skybox can be interpreted in $B$ distinct
ways, $B \rightarrow \infty$, as representations of distinct ZDs
with differing ${\bf G}$ and, therefore, differing U-indices. By
obvious analogy to the theory of Riemann surfaces in complex
analysis, each such skybox is a separate ``sheet''; as with even
such simple functions as the logarithmic, the number of such sheets
is infinite. We could then think of the infinite sequence of
skyboxes as so many cross-sections, at constant distances, of a
flashlight beam whose intensity (one over the ET's cell count)
follows Kepler's inverse square law.  Alternatively, we could ignore
the sheeting and see things another way.

Where we called fixed-$N$, ${\bf S}$ varying sequences of ETs
flip-books, we refer to fixed-${\bf S}$, $N$ varying sequences as
balloon rides:  the image is suggested by David Niven's role as
Phineas Fogg in the movie made of Jules Vernes' \textit{Around the
World in 80 Days}:  to ascend higher, David would drop a sandbag
over the side of his hot-air balloon's basket; if coming down, he
would pull a cord that released some of the balloon's steam. Each
such navigational tactic is easy to envision as a bit-shift, pushing
${\bf G}$ further to the left to cross LLs into a higher skybox,
else moving it rightward to descend. Using ${\bf S = 15}$ as the
basis of a 3-stage balloon-ride, we see how increasing $N$ from $5$
to $6$ to $7$ approaches the white-space complement of one of the
simplest (and least efficient) plane-filling fractals, the Ces\`aro
double sweep [6, p. 65].

\begin{figure}
\includegraphics[width = 1.\textwidth] {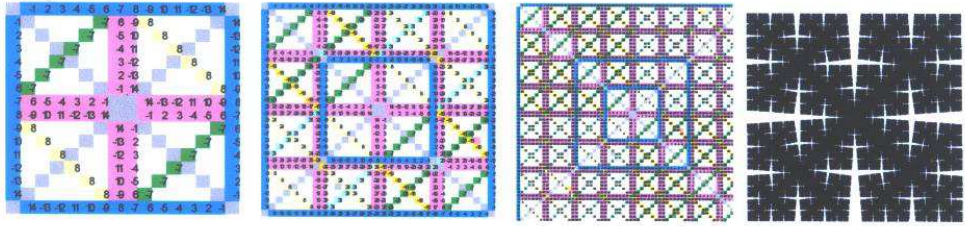}
\caption{ETs for S=15, N=5,6,7 (nested skyboxes in blue)$\cdots$ and
``fractal limit.''}
\end{figure}

The graphics were programmatically generated prior to the proving of
the theorems we're elaborating:  their empirical evidence was what
informed (indeed, demanded) the theoretical apparatus.  And we are
not quite finished with the current task the apparatus requires of
us.  We need two more theorems to finish the discussion of skybox
recursion.  For both, suppose some skybox with $B = k$, $k$ any
non-negative integer, is nested in one with $B = k + 1$.  Divide the
former along midlines to frame its four quadrants, then block out
the latter skybox into a $4 \times 4$ grid of same-sized window
panes, partitioned by the one-cell-thick borders of its own midlines
into quadrants, each of which is further subdivided by the outside
edges of the 4 one-cell-thick label lines and their extensions to
the window's frame.  These extended LLs are themselves NSLs, and
have $R, C$ values of $g$ and $g + {\bf S}$; for ${\bf S} = 15$,
they also adjoin NSLs along their outer edges whose $R, C$ values
are multiples of $8$ plus ${\bf S} \; mod \; 8$.  These pane-framing
pairs of NSLs we will henceforth refer to (as a windowmaker would)
as \textit{muntins}.  It is easy to calculate that while the inner
skybox has but one muntin each among its rows and columns, each
further nesting has $2^{B+1} - 1$.   But we are getting ahead of
ourselves, as we still have two proofs to finish. Let's begin with
Four Corners, or

\medskip

\noindent {\small Theorem 12.}  The 4 panes in the corners of the
16-paned $B = k + 1$ window are identical in contents and marks to
the analogously placed quadrants of the $B = k$ skybox.

\pagebreak

\noindent \textit{Proof.} Invoke the Zero-Padding Lemma with regard
to the U-indices, as the labels of the boxes in the corners of the
$B = k + 1$ ET are identical to those of the same-sized quadrants in
the $B = k$ ET, all labels $\geq$ the latter's $g$ only occurring in
the newly inserted region. $\;\; \blacksquare$

\medskip

\noindent \textit{Remarks.}  For $N = 6$, all filled Four Corners
cells indicate edges belonging to $3$ Box-Kites, whose edges they in
fact exhaust.  These $3$, not surprisingly, are the zero-padded
versions of the identically L-indexed trio which span the entirety
of the $N = 5$ ET. By calculations we'll see shortly, however, the
inner skybox, when considered as part of the $N = 6$ ET, has filled
cells belonging to all the other 16 Box-Kites, even though the
contents of these cells are identical to those in the $N = 5$ ET. As
$B$ increases, then, the ``sheets'' covering this same central
region must draw upon progressively more extensive networks of
interconnected Box-Kites.  As we approach the fractal limit -- and
``the Sky is the limit'' -- these networks hence become scale-free.
(Corollarily, for $N = 7$, the Four Corners' cells exhaust all the
edges of the $N = 6$ ET's 19 Box-Kites, and so on.)

Unlike a standard fractal, however, such a Sky merits the prefix
``meta'':  for each empty ET cell corresponds to a point in the
usual fractal variety; and each pair of filled ET cells, having
(r,c) of one = (c,r) of the other), correspond to diagonal-pairs in
Assessor planes, orthogonal to all other such diagonal-pairs
belonging to the other cells. Each empty ET cell, in other words,
not only corresponds to a point in the usual plane-confined fractal,
but belongs to the complement of the filled cells' infinite number
of \textit{dimensions} framing the Sky's \textit{meta-}fractal.

\medskip

We've one last thing to prove here.  The French Windows Theorem
shows us the way the cell contents of the pairs of panes contained
between the $B = k + 1$ skybox's corners are generated from those of
the analogous pairings of quadrants in the $B = k$ skybox, by adding
$g$ to L-indices.

\medskip

\noindent {\small Theorem 13}.  For each half-square array of cells
created by one or the other midline (the French windows), each cell
in the half-square parallel to that adjoining the midline (one of
the two shutters), but itself adjacent to the label-line delimiting
the former's bounds, has content equal to $g$ plus that of the cell
on the same line orthogonal to the midline, and at the same distance
from it, as \textit{it} is from the label-line.  All the empty
long-diagonal cells then map to $g$ (and are marked), or $g + {\bf
S}$ (and are unmarked).  Filled cells in extensions of the
label-lines bounding each shutter are calculated similarly, but with
reversed markings; all other cells in a shutter have the same marks
as their French-window counterparts.

\medskip

\noindent \textit{Preamble.}  Note that there can be (as we shall
see when we speak of \textit{hide/fill involution}) cells left empty
for rule-based reasons other than $P \; = \; R \veebar C \; = \; 0
\;|
 \; {\bf S}$. The shutter-based counterparts of such French-window
cells, unlike those of long-diagonal cells, remain empty.

\medskip

\noindent \textit{Proof.}  The top and left (bottom and right)
shutters are equivalent: one merely switches row for column labels.
Top/left and bottom/right shutter-sets are likewise equivalent by
the symmetry of strut-opposites.  We hence make the case for the
left shutter only.  But for the novelties posed by the initially
blank cells and the label lines (with the only real subtleties
involving markings), the proof proceeds in a manner very similar to
Theorem 11:  split into 3 cases, based on whether (1) the L-index
trip implied by the $R, C, P$ values is a Zigzag; (2) $u, v$ are
both Vents; or, (3) the edge signified by the cell content is the
emanation of same-inner-signed dyads (that is, one is a Vent, the
other a Zigzag).

\textit{Case 1:}  Assume $(u, v, w)$ a Zigzag L-trip in the French
window's contained skybox; the general product in its shutter is

\begin{center}
$v \;\; - \;\; ({\bf G} + v_{opp})$

\underline{$\;\;\;(g + u)   \;\; + \;\;({\bf G} + g +
u_{opp})\;\;\;$}
\smallskip

$- ({\bf G} + g + w_{opp}) \;\;\; + (g + w)$

\underline{$ - (g + w)  \;\;\; + ({\bf G} + g + w_{opp})$}

\smallskip
$0$

\end{center}

\smallskip

$(g + u) \cdot v = P:$  $(u,v,w) \rightarrow (g + w, v, g + u)$;
hence, $(- (g + w))$.

$(g + u) \cdot ({\bf G} + v_{opp}) = P:$  $(u, w_{opp}, v_{opp})$
$\rightarrow (g + w_{opp}, g + u, v_{opp})$ $\rightarrow ({\bf G} +
v_{opp}, g + u, {\bf G} + g + w_{opp})$; dyads' opposite inner signs
make $({\bf G} + g + w_{opp})$ positive.

$({\bf G} + g + u_{opp}) \cdot v = P:$  $(u_{opp}, w_{opp}, v)$
$\rightarrow (g + w_{opp}, g + u_{opp}, v)$ $\rightarrow ({\bf G} +
g + u_{opp}, {\bf G} + g + w_{opp}, v)$; hence, $(- ({\bf G} + g +
w_{opp}))$.

$({\bf G} + g + u_{opp}) \cdot ({\bf G} + v_{opp}) = P:$  $(v_{opp},
u_{opp}, w)$ $\rightarrow$ $(v_{opp}, g + w, g + u_{opp})$
$\rightarrow$ $({\bf G} + g + u_{opp}, g + w, {\bf G} + v_{opp})$;
dyads' opposite inner signs make $(g + w)$ positive.

\medskip

\textit{Case 2:}  The product of two Vents must have negative
edge-sign, hence negative inner sign in top dyad to lower dyad's
positive.  The shutter product thus looks like this:

\pagebreak

\begin{center}
$(u_{opp}) - ({\bf G} + u)$

\underline{$\;(g + v_{opp})   \; + \;({\bf G} + g + v)\;$}
\smallskip

$+ ({\bf G} + g + w_{opp}) \;\;\;\;\; + (g + w)$

\underline{$ - (g + w)  \;\;\;\;\; - ({\bf G} + g + w_{opp})$}

\smallskip
$0$

\end{center}

\smallskip

$(g + v_{opp}) \cdot u_{opp} = P:$  $(v_{opp}, u_{opp}, w)$
$\rightarrow (g + w, u_{opp}, g + v_{opp})$; hence, $(- (g + w))$.

$(g + v_{opp}) \cdot ({\bf G} + u) = P:$  $(v_{opp}, u, w_{opp})$
$\rightarrow (g + w_{opp}, u, g + v_{opp})$ $\rightarrow ({\bf G} +
u, {\bf G} + g + w_{opp}, g + v_{opp})$; but dyads' inner signs are
opposite, so $(- ({\bf G} + g + w_{opp}))$.

$({\bf G} + g + v) \cdot u_{opp} = P:$  $(u_{opp}, w_{opp}, v)$
$\rightarrow (u_{opp}, g + v, g + w_{opp})$ $\rightarrow (u_{opp},
{\bf G} + g + w_{opp}, {\bf G} + g + v)$; hence, $(+ ({\bf G} + g +
w_{opp}))$.

$({\bf G} + g + v) \cdot ({\bf G} + u) = P:$  $(u, v, w)$
$\rightarrow$ $(u, g + w, g + v)$ $\rightarrow$ $({\bf G} + g + v, g
+ w, {\bf G} + u)$; but dyads' inner signs are opposite, so $(+ (g +
w))$.

\medskip

\textit{Case 3:}  The product of Vent and Zigzag displays same inner
sign in both dyads; hence the following arithmetic holds:

\smallskip

\begin{center}
$(u_{opp}) + ({\bf G} + u)$

\underline{$\; (g + v)   \; + \;({\bf G} + g + v_{opp})\;$}
\smallskip

$+({\bf G} + g + w) \;\;\;\;\; +( g + w_{opp})$

\underline{$ - (g + w_{opp})  \;\;\;\;\; - ({\bf G} + g + w)$}

\smallskip
$0$

\end{center}
As with the last case in Theorem 11, we omit the term-by-term
calculations for this last case, as they should seem ``much of a
muchness'' by this point.  What is clear in all three cases is that
index values of shutter cells have same markings as their
French-window counterparts, at least for all cells which
\textit{have} markings in the contained skybox; but, in all cases,
indices are augmented by $g$.

The assignment of marks to the shutter-cells linked to blank cells
in French windows is straightforward for Type I box-kites:  since
any containing skybox must have $g > {\bf S}$, and since $g + s$ has
$g$ as its strut opposite, then the First Vizier tells us that any
$g$ must be a Vent. But then the $R, C$ indices of the cell
containing $g$ must belong to a Trefoil in such a box-kite; hence,
one is a Vent, the other a Zigzag, and $g$ must be marked.  Only if
the $R, C, P$ entry in the ET is necessarily confined to a Type II
box-kites will this not necessarily be so.  But Part II's Appendix B
made clear that Type II's are generated by \textit{excluding} $g$
from their L-indices:  recall that, in the Pathions, for all ${\bf
S}$ < 8, all and only Type II box-kites are created by placing one
of the Sedenion Zigzag L-trips on the ``Rule 0'' circle of the
PSL(2,7) triangle with $8$ in the middle (and hence excluded).  This
is a box-kite in its own right (one of the 7 ``Atlas'' box-kites
with ${\bf S} = 8$); its 3 sides are ``Rule 2'' triplets, and
generate Type II box-kites when made into zigzag L-index sets.
Conversely, all Pathion box-kites containing an '8' in an L-index
(dubbed ''strongboxes'' in Appendix B) are Type I. Whether something
peculiar might occur for large $N$ (where there might be multiple
powers of 2 playing roles in the same box-kite) is a matter of
marginal interest to present concerns, and will be left as an open
question for the present.  We merely note that, by a similar
argument, and with the same restrictions assumed, $g + {\bf S}$ must
be a Zigzag L-index, and $R, C$ either both be likewise (hence, $g +
{\bf S}$ is unmarked); or, both are Vents in a Trefoil (so $g + {\bf
S}$ must be unmarked here too).

The last detail -- reversal of label-line markings in their
$g$-augmented shutter-cell extensions -- is demonstrated as follows,
with the same caveat concerning Type II box-kites assumed to apply.
Such cells house DMZs (just swap $u$ for $g + u$ in Theorem 11's
first setup -- they form a Rule 1 trip -- and compute). The LL
extension on top has row-label $g$; that along the bottom, the
strut-opposite $g + {\bf S}$. Given trip $(u,v,w)$, the shutter-cell
index for $R, C = (g, u)$ corresponds to French-window index for $R,
C = (g, g + u)$. But $(u, g, g+u)$ is a Trefoil, since $g$ is a
Vent. So if $u$ is one too, $g + u$ isn't; hence marks are reversed
as claimed. $\;\;\blacksquare$

\section{Maximal High-Bit Singletons: (s,g)-Modularity for ${\bf 16 < S \leq 24}$}
The Whorfian Sky, having but one high bit in its strut constant, is
the simplest possible meta-fractal -- the first of an infinite
number of such infinite-dimensional zero-divisor-spanned spaces. We
can consider the general case of such singleton high-bit
recursiveness in two different, complementary ways.  First, we can
supplement the just-concluded series of theorems and proofs with a
calculational interlude, where we consider the iterative embeddings
of the Pathion Sand Mandalas in the infinite cascade of
boxes-within-boxes that a Sky oversees.  Then, we can generalize
what we saw in the Pathions to consider the phenomenology of strut
constants with singleton high-bits, which we take to be any bits
representing a power of $2 \geq 3$ if ${\bf S}$ contains low bits
(is not a multiple of $8$), else a power of $2$ strictly greater
than $3$ otherwise. Per our earlier notation, $g = {\bf G/2}$ is the
highest such singleton bit possible.  We can think of its
exponential increments -- equivalent to left-shifts in bit-string
terms -- as the side-effects of conjoint zero-padding of $N$ and
${\bf S}$. This will be our second topic in this section.

Maintaining our use of ${\bf S = 15}$ as exemplary, we have already
seen that NSLs come in quartets:  a row and column are each headed
by ${\bf S} \; mod \; g$ (henceforth, $s$) and $g$, hence $7$ and
$8$ in the Sand Mandalas.  But each recursive embedding of the
current skybox in the next creates further quartets.  Division down
the midlines to insert the indices new to the next CDP generation
induces the Sand Mandala's adjoining strut-opposite sets of $s$ and
$g$ lines (the pane-framing muntins) to be displaced to the borders
of the four corners and shutters, with the new skybox's $g$ and $g +
s$ now adjoining the old $s$ and $g$ to form new muntins, on the
right and left respectively, while $g + g/2$ (the old ${\bf G} + g)$
and its strut opposite form a third muntin along the new midlines.
Continuing this recursive nesting of skyboxes generates 1, 3, 7,
$\cdots$, $2^{B+1}-1$ row-and-column muntin pairs involving
multiples of $8$ and their supplementings by $s$, where (recalling
earlier notation) $B = 0$ for the inner skybox, and increments by
$1$ with each further nesting.  Put another way, we then have a
muntin number $\mu = (2^{N-4} - 1)$, or $4\mu$ NSL's in all.

The ET for given $N$ has $(2^{N-1}-2)$ cells in each row and column.
But NSLs divvy them up into boxes, so that each line is crossed by
$2 \mu$ others, with the 0, 2 or 4 cells in their overlap also
belonging to diagonals.  The number of cells in the overlap-free
segments of the lines, or $\omega$, is then just $4 \mu \cdot
(2^{N-1} - 2 - 2 \mu ) = 24 \mu ( \mu + 1 )$:  an integer number of
Box-Kites. For our ${\bf S = 15}$ case, the minimized line shuffling
makes this obvious:  all boxes are 6 x 6, with 2-cell-thick
boundaries (the muntins separating the panes), with $\mu$
boundaries, and $( \mu + 1)$ overlap-free cells per each row or
column, per each quartet of lines.

The contribution from diagonals, or $\delta$, is a little more
difficult, but straightforward in our case of interest:  4 sets of
$1, 2, 3, \cdots, \mu$ boxes are spanned by moving along
\textit{one} empty long diagonal before encountering the
\textit{other}, with each box contributing 6, and each overlap zone
between adjacent boxes adding 2.  Hence, $\delta = 24 \cdot (2^{N-3}
- 1) (2^{N-3} - 2)/6$ -- a formula familiar from associative-triplet
counting:  it also contributes an integer number of Box-Kites.  The
one-liner we want, then, is this:

\smallskip

\begin{center}

{\small $BK_{N,\; 8 < {\bf S} < 16} = \omega + \delta =
(2^{N-4})(2^{N-4} - 1) \; +  \; (2^{N-3} - 1)(2^{N-3} - 2)/6$}

\end{center}

\smallskip

For $N = 4, 5, 6, 7, 8, 9, 10$, this formula gives $0, 3, 19, 91,
395$, $1643, 6699$.  Add $4^{N-4}$ to each -- the immediate
side-effect of the offing of all four Rule 0 candidate trips of the
Sedenion Box-Kite exploded into the Sand Mandala that begins the
recursion -- and one gets ``d\'ej\`a vu all over again'': $1$, $7$,
$35$, $155$, $651$, $2667$, $10795$ -- the full set of Box-Kites for
${\bf S \leq 8}$.

It would be nice if such numbers showed up in unsuspected places,
having nothing to do with ZDs.  Such a candidate context does, in
fact, present itself, in Ed Pegg's regular MAA column on ``Math
Games'' focusing on ``Tournament Dice.'' [7]  He asks us, ``What
dice make a non-transitive four player game, so that if three dice
are chosen, a fourth die in the set beats all three?  How many dice
are needed for a five player non-transitive game, or more?''  The
low solution of 3 explicitly involves PSL(2,7); the next solution of
19 entails calculations that look a lot like those involved in
computing row and column headers in ETs.  No solutions to the
dice-selecting game beyond 19 are known.  The above formulae,
though, suggest the next should be 91.  Here, ZDs have no apparent
role save as dummies, like the infinity of complex dimensions in a
Fourier-series convergence problem, tossed out the window once the
solution is in hand.  Can a number-theory fractal, with
intrinsically structured cell content (something other, non-meta,
fractals lack) be of service in this case -- and, if not in this
particular problem, in others like it?

Now let's consider the more general situation, where the singleton
high-bit can be progressively left-shifted. Reverting to the use of
the simplest case as exemplary, use ${\bf S = g + 1 = 9}$ in the
Pathions, then do tandem left-shifts to produce this sequence: $N =
6, \; {\bf S = g + 1}$ ${\bf = 17}$; $N = 7$, ${\bf S = }$ ${\bf g +
1 = 33};\; \cdots;\; N = K$, ${\bf S = g + 1} = 2^{K-2} + 1$. A
simple rule governs these ratchetings:  in all cases, the number of
filled cells = $6 \cdot (2^{N-1} - 4)$, since there are two sets of
parallel sides which are filled but for long-diagonal intersections,
and two sets of $g$ and $1$ entries distributed one per row along
orthogonals to the empty long diagonals.  Hence, for the series just
given, we have cell counts of $72,\; 168,\; \cdots,\; 6 \cdot (2^{N
- 1} - 4)$ for $BK_{N,\;S} = 3,\;7,\; \cdots,\;2^{N - 3} - 1$, for
$g < {\bf S} < g + 8 = {\bf G}$ in the Pathions, and all $g < {\bf
S} \leq g + 8$ in the Chingons, $2^{7}$-ions, and general
$2^{N}$-ions, in that order.

Algorithmically, the situation is just as easy to see:  the
splitting of dyads, sending U- and L- indices to strut-opposite
Assessors, while incorporating the ${\bf S}$ and ${\bf G}$ of the
current CDP generation as strut-opposites in the next, continues.
For ${\bf S = 17}$ in the Chingons, there are now $2^{N-3}-1 = 7$,
not $3$, Box-Kites sharing the new $g = 16$ (at B) and ${\bf S}
\;mod \;g = 1$ (at E) in our running example. The U- indices of the
Sand Mandala Assessors for ${\bf S = g + 1 = 9}$ are now L-indices,
and so on:  every integer $< G$ and $\neq {\bf S}$ gets to be an
L-index of one of the $30 (= 2^{N-1} - 2)$ Assessors, as $16$ and
${\bf S} \;mod \; g = 1$ appear in each of the $7$ Box-Kites, with
each other eligible integer appearing once only in one of the $7
\cdot 4 = 28$ available L-index slots.

As an aside, in all 7 cases, writing the smallest Zigzag L-index at
$a$ mandates all the Trefoil trips be ``precessed'' -- a phenomenon
also observed in the ${\bf S = 8}$ Pathion case, as tabulated on p.
14 of [8]. For Zigzag L-index set $(2, 16, 18)$, for instance,
$(a,d,e)$ $=$ $(2,3,1)$ instead of $(1,2,3)$; $(f,c,e)$ $=$
$(19,18,1)$ not $(1,19,18)$; and $(f,d,b)$ $=$ $(19,3,16)$. But
otherwise, there are no surprises: for $N=7$, there are $(2^{7 - 3}
- 1) = 15$ Box-Kites, with all $62 ( = 2^{N-1} - 2)$ available cells
in the rows and columns linked to labels $g$ and ${\bf S}\; mod \;
g$ being filled, and so on.

Note that this formulation obtains for any and all ${\bf S > 8}$
where the maximum high-bit (that is, $g$) is included in its
bitstring: for, with $g$ at B and ${\bf S} \;mod \; g$ at E,
whichever ${\bf R, C}$ label is not one of these suffices to
completely determine the remaining Assessor L-indices, so that no
other  bits in ${\bf S}$ play a role in determining any of them.
Meanwhile, cell \textit{contents} ${\bf P}$ containing either $g$ or
${\bf S} \; mod \; g$, but created by XORing of row and column
labels equal to neither, are arrayed in off-diagonal pairs, forming
disjoint sets parallel or perpendicular to the two empty ones.  If
we write ${\bf S} \;mod \; g$ with a lower-case $s$, then we could
call the rule in play here \textit{(s,g)-modularity}. Using the
vertical pipe for logical or, and recalling the special handling
required by the 8-bit when ${\bf S}$ is a multiple of 8 (which we
signify with the asterisk suffixed to ``mod''), we can shorthand its
workings this way:

\medskip

\noindent {\small Theorem 14}.  For a $2^{N}$-ion inner skybox whose
strut constant ${\bf S}$ has a singleton high-bit which is maximal
(that is, equal to $g = {\bf G/2} = 2^{N-2}$), the recipe for its
filled cells can be condensed thus:

\smallskip

\begin{center}
 ${\bf R\; |\; C |\; P} = g\; |\; {\bf
S} \; mod^{*} \; g$

\end{center}

Under recursion, the recipe needs to be modified so as to include
not just the inner-skybox $g$ and ${\bf S} \; mod^{*} \; g$
(henceforth, simply lowercase $s$), but all integer multiples $k$ of
$g$ less than the ${\bf G}$ of the outermost skybox, plus their
strut opposites $k \cdot g + s$.

\medskip
\noindent \textit{Proof}.  The theorem merely boils down the
computational arguments of prior paragraphs in this section, then
applies the last section's recursive procedures to them.  The first
claim of the proof is identical to what we've already seen for Sand
Mandalas, with zero-padding injected into the argument.  The second
claim merely assumes the area quadrupling based on midline
splitting, with the side-effects already discussed.  No formal
proof, then, is called for beyond these points. $\;\;\blacksquare$

\medskip

\noindent \textit{Remarks}.  Using the computations from two
paragraphs prior to the theorem's statement, we can readily
calculate the box-kite count for any skybox, no matter how deeply
nested:  recall the formula $6 \cdot (2^{N - 1} - 4)$ for
$BK_{N,\;S} = 2^{N - 3} - 1$.  It then becomes a straightforward
matter to calculate, as well, the limiting ratio of this count to
the maximal full count possible for the ET as $N \rightarrow
\infty$, with each cell approaching a point in a standard 2-D
fractal.  Hence, for any ${\bf S}$ with a singleton high-bit in
evidence, there exists a Sky containing all recursive redoublings of
its inner skybox, and computations like those just considered can
further be used to specify fractal dimensions and the like. (Such
computations, however, will not concern us.)  Finally, recall that,
by spectrographic equivalence, all such computations will lead to
the same results for each ${\bf S}$ value in the same spectral band
or octave.

\section{Hide/Fill Involution:  Further-Right High-Bits with ${\bf 24 < S < 32}$.}
Recall that, in the Sand Mandala flip-book, each increment of ${\bf
S}$ moved the two sets of orthogonal parallel lines one cell closer
toward their opposite numbers:  while ${\bf S = 9}$ had two
filled-in rows and columns forming a square missing its corners, the
progression culminating in ${\bf S = 15}$ showed a cross-hairs
configuration:  the parallel lines of cells now abutted each other
in 2-ply horizontal and vertical arrays.  The same basic progression
is on display in the Chingons, starting with ${\bf S = 17}$.  But
now the number of strut-opposite cell pairs in each row and column
is 15, not 7, so the cross-hairs pattern can't arise until ${\bf S =
31}$.  Yet it never arises in quite the manner expected, as
something quite singular transpires just after flipping past the ET
in the middle, for ${\bf S = 24}$.  Here, rows and columns labeled
$8$ and $16$ constrain a square of empty cells in the center
$\cdots$ quickly followed by an ET which seems to continue the
expected trajectory -- except that almost all the non-long-diagonal
cells left empty in its predecessor ETs are now inexplicably filled.
More, there is a method to the ``almost all'' as well:  for we now
see not 2, but 4 rows and columns, all being blanked out while those
labeled with $g$ and ${\bf S} \; mod \; g$ are being filled in.

This is an inevitable side effect of a second high-bit in ${\bf S}$:
we call this phenomenon, first appearing in the Chingons,
\textit{hide/fill involution}.  There are 4, not 2, line-pairs,
because ${\bf S}$ and ${\bf G}$, modulo a lower power of 2 (because
devolving upon a prior CDP generation's $g$), offer twice the
possibilities: for ${\bf S = 25}$, ${\bf S} \; mod \; 16$ is now
$9$, but ${\bf S} \; mod\; 8$ can result in either $1$ or $17$ as
well -- with correlated \textit{multiples} of $8$ ($8$ proper, and
$24$) defining the other two pairings.  All cells with ${\bf R \; |
C \; | \; P}$ equal to one of these 4 values, but for the handful
already set to ``on'' by the first high-bit, will now be set to
``off,'' while all other non-long-diagonal cells set to ``off'' in
the Pathion Sand Mandalas are suddenly ``on.''  What results for
each Chingon ET with $24 < {\bf S} < 32$ is an ensemble comprised of
$23$ Box-Kites. (For the flip-book, see Slides 40 -- 54.)  Why does
this happen? The logic is as straightforward as the effect can seem
mysterious, and is akin, for good reason, to the involutory effect
on trip orientation induced by Rule 2 addings of ${\bf G}$ to 2 of
the trip's 3 indices.

In order to grasp it, we need only to consider another pair of
abstract calculation setups, of the sort we've seen already many
times.  The first is the core of the Two-Bit Theorem, which we state
and prove as follows:

\medskip

\noindent {\small Theorem 15}.  $2^{N}$-ion dyads making DMZs before
augmenting ${\bf S}$ with a new high-bit no longer do so after the
fact.

\medskip

\noindent \textit{Proof}. Suppose the high-bit in the bitstring
representation of ${\bf S}$ is $2^{K},\;K < (N-1)$. Suppose further
that, for some L-index trip $(u,v,w)$, the Assessors $U$ and $V$ are
DMZ's, with their dyads having same inner signs. (This last
assumption is strictly to ease calculations, and not substantive: we
could, as earlier, use one or more binary variables of the $sg$ type
to cover all cases explicitly, including Type I vs. Type II
box-kites.  To keep things simple, we assume Type I in what
follows.) We then have $(u + u \cdot X)(v + v \cdot X) = (u + U)(v +
V) = 0$. But now suppose, without changing $N$, we add a bit
somewhere further to the left to ${\bf S}$, so that ${\bf S} <
(2^{K} = L) < {\bf G}$. The augmented strut constant now equals
${\bf S_{L}} = {\bf S + L}$. One of our L-indices, say $v$, belongs
to a Vent Assessor thanks to the assumed inner signing; hence, by
Rule 2 and the Third Vizier, $(V,v,X) \rightarrow (X + L, v, V +
L)$.  Its DMZ partner $u$, meanwhile, must thereby be a Zigzag
L-index, which means $(u,U,X) \rightarrow (u, X + L, U + L)$. We
claim the truth of the following arithmetic:

\smallskip

\begin{center}

$v \; + \; (V + L)$ \\
\underline{$ \;\;\; u \; + \; (U + L)\;\;\;$} \\
$+(W \;+ \; L) \; + w$ \\
\underline{$+ \; w \;\; - (W +  L)$} \\
NOT ZERO (+w's don't cancel) \\

\end{center}

\smallskip

The left bottom product is given. The product to its right is
derived as follows:  since $u$ is a Zigzag L-index, the Trefoil
U-trip $(u,V,W)$ has the same orientation as $(u,v,w)$, so that Rule
2 $\rightarrow (u, W+L, V+L)$, implying the negative result shown.
The left product on the top line, though, has terms derived from a
Trefoil U-trip lacking a Zigzag L-index, so that only after Rule 2
reversal are the letters arrayed in Zigzag L-trip order:  $(U + L,
v, W + L)$.  Ergo, $+ (W+L)$.  Similarly for the top right:  Rule 2
reversal ``straightens out'' the Trefoil U-trip, to give $(U + L, V
+ L, w)$; therefore, $(+ w)$ results.  If we explicitly covered
further cases by using an $sg$ variable, we would be faced with a
Theorem 2 situation:  one or the other product pair cancels, but not
both. $\;\;\blacksquare$

\medskip

\noindent \textit{Remark.}  The prototype for the phenomenon this
theorem covers is the ``explosion'' of a Sedenion box-kite into a
trio of interconnected ones in a Pathion sand mandala, with the
${\bf S}$ of the latter = the ${\bf X}$ of the former.  As part of
this process, 4 of the expected 7 are ``hidden'' box-kites (HBKs),
with no DMZs along their edges.  These have zigzag L-trips which are
precisely the L-trips of the 4 Sedenion Sails. Here, an empirical
observation which will spur more formal investigations in a sequel
study:  for the 3 HBKs based on trefoil L-trips, exactly 1 strut has
reversed orientation (a different one in each of them), with the
orientation of the triangular side whose midpoint it ends in also
being reversed.  For the HBK based on the zigzag L-trip, all 3
struts are reversed, so that the flow along the sides is exactly the
reverse of that shown in the ``Rule 0'' circle.  (Hence, all
possible flow patterns along struts are covered, with only those
entailing 0 or 2 reversals corresponding to functional box-kites:
our Type I and Type II designations.)  It is not hard to show that
this zigzag-based HBK has another surprising property:  the 8 units
defined by its own zigzag's Assessors plus ${\bf X}$ and the real
unit form a ZD-free copy of the Octonions.  This is also true when
the analogous Type II situation is explored, albeit for a slightly
different reason:  in the former case, all 3 Catamaran ``twistings''
take the zigzag edges to other HBKs; in the latter, though, the pair
of Assessors in some other Type II box-kite reached by ``twisting''
-- $(a,B)$ and $(A,b)$, say, if the edge be that joining Assessors A
and B, with strut-constant $c_{opp} = d$ -- are \textit{strut
opposites}, and hence also bereft of ZDs.  The general picture seems
to mirror this concrete case, and will be studied in ``Voyage by
Catamaran'' with this expectation:  the bit-twiddling logic that
generates meta-fractal ``Skies'' also underwrites a means for
jumping between ZD-free Octonion clones in an infinite number of
HBKs housed in a Sky.  Given recent interest in pure ``E8'' models
giving a privileged place to the basis of zero-divisor theory,
namely ``G2'' projections (viz., A. Garrett Lisi's ``An
Exceptionally Simple Theory of Everything''); a parallel vogue for
many-worlds approaches; and, the well-known correspondence between
8-D closest-packing patterns, the loop of the 240 unit Octonions
which Coxeter discovered, and E8 algebras -- given all this,
tracking the logic of the links across such Octonionic ``brambles''
might prove of great interest to many researchers.

\medskip

Now, we still haven't explained the flipside of this off-switch
effect, to which prior CDP generation Box-Kites -- appropriately
zero-padded to become Box-Kites in the current generation until the
new high-bit is added to the strut-constant -- are subjected.  How
is it that previously empty cells \textit{not} associated with the
second high-bit's blanked-out R, C, P values are now \textit{full}?
The answer is simple, and is framed in the Hat-Trick Theorem this
way.

\medskip

\noindent {\small Theorem 16}.  Cells in an ET which represent DMZ
edges of some $2^{N}$-ion Box-Kites for some fixed ${\bf S}$, and
which are offed in turn upon augmenting of ${\bf S}$ by a new
leftmost bit, are turned on once more if ${\bf S}$ is augmented by
yet another new leftmost bit.

\medskip

\noindent \textit{Proof}.  We begin an induction based upon the
simplest case (which the Chingons are the first $2^{N}$-ions to
provide): consider Box-Kites with ${\bf S \leq 8}$. If a high-bit be
appended to ${\bf S}$, then the associated Box-Kites are offed.
However, if \textit{another} high-bit be affixed, these dormant
Box-Kites are re-awakened -- the second half of \textit{hide/fill
involution}.  We simply assume an L-index set $(u,v,w)$ underwriting
a Sail in the ET for the pre-augmented ${\bf S}$, with Assessors
$(u, U)$ and $(v, V)$.  Then, we introduce a more leftified bit
$2^{Q} = M$, where pre-augmented ${\bf S} < L < M < {\bf G}$, then
compute the term-by term products of $(u + (U + L + M))$ and $(v +
sg \cdot (V + L + M))$, using the usual methods. And as these
methods tell us that two applications of Rule 2 have the same effect
as none in such a setup, we have no more to prove.
$\;\;\blacksquare$

\medskip

\noindent \textit{Corollary}.  The induction just invoked makes it
clear that strut constants equal to multiples of $8$ not powers of
$2$ are included in the same spectral band as all other integers
larger than the prior multiple.  The promissory note issued in the
second paragraph of Part II's concluding section, on 64-D
Spectrography, can now be deemed redeemed.

\medskip

In the Chingons, high-bits $L$ and $M$ are necessarily adjacent in
the bitstring for ${\bf S < G = 32}$; but in the general $2^{N}$-ion
case, $N$ large, zero-padding guarantees that things will work in
just the same manner, with only one difference:  the recursive
creation of ``harmonics'' of relatively small-$g$ $(s,g)$-modular
${\bf R, C, P}$ values will propagate to further levels, thereby
effecting overall Box-Kite counts.

In general terms, we have echoes of the formula given for
$(s,g)$-modular calculations, but with this signal difference: there
will be \textit{one} such rule for \textit{each} high-bit $2^{H}$ in
${\bf S}$, where residues of ${\bf S}$ modulo $2^{H}$ will generate
their own near-solid lines of rows and columns, be they hidden or
filled.  Likewise for multiples of $2^{H} < {\bf G}$ which are not
covered by prior rules, and multiples of $2^H$ supplemented by the
bit-specific residue (regardless of whether $2^{H}$ itself is
available for treatment by this bit-specific rule). In the simplest,
no-zero-padding instances, all even multiples are excluded, as they
will have occurred already in prior rules for higher bits, and fills
or hides, once fixed by a higher bit's rule, cannot be overridden.

Cases with some zero-padding are not so simple.  Consider this
two-bit instance, ${\bf S = 73}, N = 8$:  the fill-bit is 64, the
hide-bit is just 8, so that only 9 and 64 generate NSLs of filled
values; all other multiples of 8, and their supplementing by 1
(including 65) are NSLs of hidden values. Now look at a variation on
this example, with the single high-bit of zero-padding removed --
i.~e., ${\bf S = 41}, N = 8$. Here, the fill-bit is 32, and its
multiples 64 and 96, as well as their supplements by ${\bf S} \;
modulo \; 32 \; = \; 9$, or 9 and 73 and 105, label NSLs of filled
values; but all other multiples of 8, plus all multiples of 8
supplemented by 1 not equal to 9 or 73 or 105, label NSLs of hidden
values.  Cases with multiple fill and hide bits, with or without
additional zero-padding, are obviously even more complicated to
handle explicitly on a case-by-case basis, but the logic framing the
rules remain simple; hence, even such messy cases are
programmatically easy to handle.

Hide/fill involution means, then, that the first, third, and any
further odd-numbered high-bits (counting from the left) will
generate ``fill'' rules, whereas all the even-numbered high-bits
generate ``hide'' rules -- with all cells not touched by a rule
being either hidden (if the total number of high-bits $B$ is odd) or
filled ($B$ is even).

Two further examples should make the workings of this protocol more
clear. First, the Chingon test case of ${\bf S = 25}$:  for $({\bf R
\; | \; C \; | \; P} \; = \; 9 \;| \; 16)$, all the ET cells are
filled; however, for $({\bf R \; | \; C \; | \; P} \; = 1 \;|\; 8
\;|\; 17 \;|\; 24)$, ET cells not already filled by the first rule
(and, as visual inspection of Slide 48 indicates, there are only 8
cells in the entire 840-cell ET already filled by the prior rule
which the current rule would like to operate on) are hidden from
view. Because the 16- and 8- bits are the only high-bits, the count
of same is even, meaning all remaining ET cells not covered by these
2 rules are filled.

We get 23 for Box-Kite count as follows.  First, the 16-bit rule
gives us 7 Box-Kites, per earlier arguments; the 8-bit rule, which
gives 3 filled Box-Kites in the Pathions, recursively propagates to
cover 19 hidden Box-Kites in the Chingons, according to the formula
produced last section. But hide/fill involution says that, of the 35
maximum possible Box-Kites in a Chingon ET, $35 - 19 = 16$ are now
made visible.  As none of these have the Pathion ${\bf G = 16}$ as
an L-index, and all the 7 Box-Kites from the 16-bit rule
\textit{do}, we therefore have a grand total of $7 + 16 = 23$
Box-Kites in the ${\bf S =25}$ ET, as claimed (and as cell-counting
on the cited Slide will corroborate).

The concluding Slides 76--78 present a trio of color-coded
``histological slices'' of the hiding and filling sequence
(beginning with the blanking of the long diagonals) for the simplest
3-high-bit case, $N = 7, {\bf S = 57}$.  Here, the first fill rule
works on 25 and 32; the first hide rule, on 9, 16, 41, and 48; the
second fill rule, on 1, 8, 17, 24, 33, 40, 49, and 56; and the rest
of the cells, since the count of high-bits is odd, are left blank.

We do not give an explicit algorithmic method here, however, for
computing the number of Box-Kites contained in this 3,720-cell ET.
Such recursiveness is best handled programmatically, rather than by
cranking out an explicit (hence, long and tedious) formula, meant
for working out by a time-consuming hand calculation.  What we can
do, instead, is conclude with a brief finale, embodying all our
results in the simple ``recipe theory'' promised originally, and
offer some reflections on future directions.

\section{Fundamental Theorem of Zero-Divisor Algebra}

All of the prior arguments constitute steps sufficient to
demonstrate the Fundamental Theorem of Zero-Divisor Algebra.  Like
the role played by its Gaussian predecessor in the legitimizing of
another ``new kind of [complex] number theory,'' its simultaneous
simplicity and generality open out on extensive new vistas at once
alien and inviting.  The Theorem proper can be subdivided into a
Proposition concerning all integers, and a ``Recipe Theory''
pragmatics for preparing and ``cooking'' the meta-fractal entities
whose existence the proposition asserts, but cannot tell us how to
construct.

\medskip

\noindent \textit{Proposition:}  Any integer $K > 8$ not a power of
$2$ can uniquely be associated with a Strut Constant ${\bf S}$ of ZD
ensembles, whose inner skybox resides in the $2^{N}$-ions with
$2^{N-2} < K < 2^{N-1}$.  The bitstring representation of ${\bf S}$
completely determines an infinite-dimensional analog of a standard
plane-confined fractal, with each of the latter's points associated
with an empty cell in the infinite Emanation Table, with all
non-empty cells comprised wholly of mutually orthogonal primitive
zero-divisors, one line of same per cell.

\medskip

\noindent \textit{Preparation:}  Prepare each suitable ${\bf S}$ by
producing its bitstring representation, then determining the number
of high-bits it contains:  if ${\bf S}$ is a multiple of 8,
right-shift 4 times; otherwise, right-shift 3 times.  Then count the
number $B$ of 1's in the shortened bitstring that results.  For this
set \verb|{B}| of $B$ elements, construct two same-sized arrays,
whose indices range from $1$ to $B$:  the array \verb|{i}| which
indexes the left-to-right counting order of the elements of
\verb|{B}|; and, the array \verb|{P}| which indexes the powers of
$2$ of the same element in the same left-to-right order.  (Example:
if $K = 613$, the inner skybox is contained in the $2^{11}$-ions; as
the number is not a multiple of $8$, the bistring representation
$1001100101$ is right-shifted thrice to yield the substring of
high-bits $1001100$; $B = 3$, and for $1 \; \leq \; i \; \leq \; 3$,
$P_{1} = 9, \; P_{2} = 6; P_{3} = 5$.)

\medskip

\noindent \textit{Cookbook Instructions:}  \begin{description}

\item \verb|[0]| ~ For a given strut-constant ${\bf S}$, compute
the high-bit count $B$ and bitstring arrays \verb|{i}| and
\verb|{P}|, per preparation instructions.

\item \verb|[1]|~ Create a square spreadsheet-cell array, of edge-length $2^{I}$,
where $I \geq {\bf G/2} = g$ of the inner skybox for ${\bf S}$, with
the Sky as the limit when $I \rightarrow \infty$.

\item \verb|[2]| ~ Fill in the labels along all four edges, with those running
along the right (bottom) borders identical to those running along
the left (top), except in reversed left-right (top-bottom) order.
Refer to those along the top as column numbers $C$, and those along
the left edge, as row numbers $R$, setting candidate contents of any
cell (r,c) to $R \veebar C = P$.

\item \verb|[3]| ~ Paint all cells along the long diagonals of the
spreadsheet just constructed a color indicating BLANK, so that all
cells with $R = C$ (running down from upper left corner) else $R
\veebar C = {\bf S}$ (running down from upper right) have their
$P$-values hidden.

\item \verb|[4]| ~ For $1 \; \leq \; i \; \leq \; B$, consider for painting only
those cells in the spreadsheet created in \verb|[1]| with $R \; | \;
C \; | \; P \; = \; m \cdot 2^{\gamma} \; | \; m \cdot 2^{\gamma} \;
+ \; \sigma$, where $\gamma = P_{i}, \sigma = {\bf S} \; mod* \;
2^{\gamma}$, and $m$ is any integer $\geq 0$ (with $m = 0$ only
producing a legitimate candidate for the right-hand's second option,
as an XOR of $0$ indicates a long-diagonal cell).

\item \verb|[5]| ~ If a candidate cell has already been painted
by a prior application of these instructions to a prior value of
$i$, leave it as is.  Otherwise, paint it with $R \veebar C$ if $i =
$ odd, else paint it BLANK.

\item \verb|[6]| ~ Loop to \verb|[4]| after incrementing $i$.  If $i
< B$, proceed until this step, then reloop, reincrement, and retest
for $i = B$.  When this last condition is met, proceed to the next
step.

\item \verb|[7]| ~ If $B$ is odd, paint all cells not already
painted, BLANK; for $B$ even, paint them with $R \veebar C$.

\end{description}

\medskip

In these pseudocode instructions, no attention is given to edge-mark
generation, performance optimization, or other embellishments.
Recursive expansion beyond the chosen limits of the $2^{N}$-ion
starting point is also not addressed. (Just keep all painted cells
as is, then redouble until the expanded size desired is attained;
compute appropriate insertions to the label lines, then paint all
new cells according to the same recipe.) What should be clear,
though, is any optimization cannot fail to be qualitatively more
efficient than the code in the appendix to [9], which computes on a
cell-by-cell basis. For ${\bf S} > 8, \; N > 4$, we've reached the
onramp to the Metafractal Superhighway:  new kinds of efficiency,
synergy, connectedness, and so on, would seem to more than
compensate for the increase in dimension.

It is well-known that Chaotic attractors are built up from fractals;
hence, our results make it quite thinkable to consider Chaos Theory
from the vantage of pure Number $\cdots$ and hence the switch from
one mode of Chaos to another as a bitstring-driven -- or, put
differently, a cellular automaton-type -- process, of Wolfram's
Class 4 complexity.  Such switching is of the utmost importance in
coming to terms with the most complex finite systems known: human
brains.  The late Francisco Varela, both a leading visionary in
neurological research and its computer modeling, and a long-time
follower of Madhyamika Buddhism who'd collaborated with the Dalai
Lama in his ``Tibetan Buddhists talk with brain scientists''
dialogues [10], pointed to just the sorts of problems being
addressed here as the next frontier.  In a review essay he
co-authored in 2001 just before his death [11, p. 237], we read
these concluding thoughts on the theme of what lies ``Beyond
Synchrony'' in the brain's workings:

\begin{quote}
The transient nature of coherence is central to the entire idea of
large-scale synchrony, as it underscores the fact that the system
does not behave dynamically as having stable attractors [e.g.,
Chaos], but rather metastable patterns -- a succession of
self-limiting recurrent patterns.  In the brain, there is no
``settling down'' but an ongoing change marked only by transient
coordination among populations, as the attractor itself changes
owing to activity-dependent changes and modulations of synaptic
connections.
\end{quote}

Varela and Jean Petitot (whose work was the focus of the intermezzo
concluding Part I, in which semiotically inspired context the Three
Viziers were introduced) were long-time collaborators, as evidenced
in the last volume on \textit{Naturalizing Phenomenology} [12] which
they co-edited. It is only natural then to re-inscribe the theme of
mathematizing semiotics into the current context: Petitot offers
separate studies, at the ``atomic'' level where Greimas' ``Semiotic
Square'' resides; and at the large-scale and architectural, where
one must place L\'evi-Strauss's ``Canonical Law of Myths.''  But the
pressing problem is finding a smooth approach that lets one slide
the same modeling methodology from the one scale to the other:  a
fractal-based ``scale-free network'' approach, in other words.  What
makes this distinct from the problem we just saw Varela consider is
the focus on the structure, rather than dynamics, of transient
coherence -- a focus, then, in the last analysis, on a
characterization of \textit{database architecture} that can at once
accommodate meta-chaotic transiency and structural linguists'
cascades of ``double articulations.''

Starting at least with C. S. Peirce over a century ago, and
receiving more recent elaboration in the hands of J. M. Dunn and the
research into the ``Semantic Web'' devolving from his work, data
structures which include metadata at the same level as the data
proper have led to a focus on ``triadic logic,'' as perhaps best
exemplified in the recent work of Edward L. Robertson. [13]  His
exploration of a natural triadic-to-triadic query language deriving
from Datalog, which he calls Trilog, is not (unlike our Skies)
intrinsically recursive.  But his analysis depends upon recursive
arguments built atop it, and his key constructs are strongly
resonant with our own (explicitly recursive) ones.  We focus on just
a few to make the point, with the aim of provoking interest in
fusing approaches, rather than in proving any particular results.

The still-standard technology of relational databases based on SQL
statements (most broadly marketed under the Oracle label) was itself
derived from Peirce's triadic thinking:  the creator of the
relational formalism, Edgar F. ``Ted'' Codd, was a PhD student of
Peirce editor and scholar Arthur W. Burks.  Codd's triadic
``relations,'' as Robertson notes (and as Peirce first recognized,
he tells us, in 1885), are ``the minimal, and thus most uniform''
representations ``where metadata, that is data about data, is
treated uniformly with regular data.''  In Codd's hands (and in
those of his market-oriented imitators in the SQL arena), metadata
was ``relegated to an essentially syntactic role'' [13, p. 1] -- a
role quite appropriate to the applications and technological
limitations of the 1970's, but inadequate for the huge and/or highly
dynamic schemata that are increasingly proving critical in
bioinformatics, satellite data interpretation, Google server-farm
harvesting, and so on. As Robertson sums up the situation motivating
his own work,

\begin{quote}
Heterogeneous situations, where diverse schemata represent
semantically similar data, illustrate the problems which arise when
one person's semantics is another's syntax -- the physical ``data
dependence'' that relational technology was designed to avoid has
been replaced by a structural data dependence.  Hence we see the
need to [use] a simple, uniform relational representation where the
data/metadata distinction is not frozen in syntax. [13, pp. 1-2]
\end{quote}

 As in relational database theory and practice, the
forming and exploiting of inner and outer \textit{joins} between
variously keyed tables of data is seminal to Robertson's approach as
well as Codd's. And while the RDF formalism of the Semantic Web (the
representational mechanism for describing structures as well as
contents of web artifacts on the World Wide Web) is likewise
explicitly triadic, there has, to date, been no formal mechanism put
in place for manipulating information in RDF format.  Hence, ``there
is no natural way to restrict output of these mechanisms to triples,
except by fiat'' [13, p. 4], much less any sophisticated rule-based
apparatus like Codd's ``normal forms'' for querying and tabulating
such data.  It is no surprise, then, that Robertson's ``fundamental
operation on triadic relations is a particular three-way join which
takes explicit advantage of the triadic structure of its operands.''
This \textit{triadic join}, meanwhile, ``results in another triadic
relation, thus providing the closure required of an algebra.'' [13,
p. 6]

Parsing Robertson's compact symbolic expressions into something
close to standard English, the trijoin of three triadic relations R,
S, T is defined as some $(a, b, c)$ selected from the universe of
possibilities $(x, y, z)$, such that $(a, x, z) \in R$, $(x, b, y)
\in S$, and $(z, y, c) \in T$.  This relation, he argues, is the
most fundamental of all the operators he defines.  When supplemented
with a few constant relations (analogs of Tarski's  ``infinite
constants'' embodied in the four binary relations of universality of
all pairs, identity of all equal pairs, diversity of all unequal
pairs, and the empty set), it can express all the standard monotonic
operators (thereby excluding, among his primitives, only the
relative complement).

How does this compare with our ZD setup, and the workings of Skies?
For one thing, Infinite constants, of a type akin to Tarski's, are
embodied in the fact that any full meta-fractal requires the use of
an infinite ${\bf G}$, which sits atop an endless cascade of
singleton leftmost bits, determining for any given ${\bf S}$ an
indefinite tower of ZDs.  One of the core operators massaging
Robertson's triads is the \textit{flip}, which fixes one component
of a relation while interchanging the other two $\cdots$ but our
Rule 2 is just the recursive analog of this, allowing one to move up
and down towers of values with great flexibilty (allowing, as well,
on and off switching effecting whole ensembles).  The integer triads
upon which our entire apparatus depends are a gift of nature, not
dictated ``by fiat,'' and give us a natural basis for generating and
tracking unique IDs with which to ``tag'' and ``unpack'' data (with
``storage'' provided free of charge by the empty spaces of our
meta-fractals:  the ``atoms'' of Semiotic Squares have four
long-diagonal slots each, one per each of the ``controls'' Petitot's
Catastrophe Theory reading calls for, and so on.)

Finally, consider two dual constructions that are the core of our
own triadic number theory: if the $(a, b, c)$ of last paragraph, for
instance, be taken as a Zigzag's L-index set, then the other trio of
triples correlates quite exactly with the Zigzag U-trips. And this
3-to-1 relation, recall, exactly parallels that between the 3
Trefoil, and 1 Zigzag, Sails defining a Box-Kite, with this very
parallel forming the support for the recursion that ultimately lifts
us up into a Sky.  We can indeed make this comparison to Robinson's
formalism exceedingly explicit:  if his X, Y, Z be considered the
angular nodes of PSL(2,7) situated at the 12 o'clock apex and the
right and left corners respectively, then his $(a, b, c)$ correspond
exactly to our own Rule 0 trip's same-lettered indices!

Here, we would point out that these two threads of reflection -- on
underwriting Chaos with cellular-automaton-tied Number Theory, and
designing new kinds of database architectures -- are hardly
unrelated.  It should be recalled that two years prior to his
revolutionary 1970 paper on relational databases [14], Codd
published a pioneering book on cellular automata [15].  It is also
worth noting that one of the earliest technologies to be spawned by
fractals arose in the arena of data compression of images, as
epitomized in the work of Michael Barnsley and his Iterative Systems
company. The immediate focus of the author's own commercial efforts
is on fusing meta-fractal mathematics with the context-sensitive
adaptive-parsing ``Meta-S'' technology of business associate Quinn
Tyler Jackson. [16] And as that focus, tautologically, is not
mathematical \textit{per se}, we pass it by and leave it, like so
many other themes just touched on here, for later work.

\pagebreak

\section*{References}
\begin{description}

\item \verb|[1]| Robert P. C. de Marrais, ``Placeholder
Substructures I:  The Road From NKS to Scale-Free Networks is Paved
with Zero Divisors,'' \textit{Complex Systems}, 17 (2007), 125-142;
arXiv:math.RA/0703745

\item \verb|[2]| Robert P. C. de Marrais, ``Placeholder
Substructures II:  Meta-Fractals, Made of Box-Kites, Fill
Infinite-Dimensional Skies,'' arXiv:0704.0026 [math.RA]

\item \verb|[3]| Robert P. C. de Marrais, ``The 42 Assessors and
the Box-Kites They Fly,'' arXiv:math.GM/0011260

\item \verb|[4]| ~ Robert P. C. de Marrais, ``The Marriage of Nothing
and All:  Zero-Divisor Box-Kites in a `TOE' Sky,'' in Proceedings of
the $26^{\textrm{th}}$ International Colloquium on Group Theoretical
Methods in Physics, The Graduate Center of the City University of
New York, June 26-30, 2006, forthcoming from Springer--Verlag.

\item \verb|[5]| Robert P. C. de Marrais, ``Placeholder
Substructures:  The Road from NKS to Small-World, Scale-Free
Networks Is Paved with Zero-Divisors,'' http:// \newline
wolframscience.com/conference/2006/
presentations/materials/demarrais.ppt  (Note:  the author's surname
is listed under ``M,'' not ``D.'')

\item \verb|[6]|~ Benoit Mandelbrot, \textit{The Fractal Geometry of Nature} (W.
H. Freeman and Company, San Francisco, 1983)

\item \verb|[7]| Ed Pegg, Jr., ``Tournament Dice,'' \textit{Math Games} column for July 11, 2005, on the MAA
website at http://www.maa.org/editorial/   mathgames/mathgames
\verb|_07_11_05|.html

\item \verb|[8]| Robert P. C. de Marrais, ``The `Something From
Nothing' Insertion Point,''
http://www.wolframscience.com/conference/2004/ presentations/
\newline
materials/rdemarrais.pdf

\item \verb|[9]| Robert P. C. de Marrais, ``Presto! Digitization,'' arXiv:math.RA/0603281

\item \verb|[10]| Francisco Varela, editor, \textit{Sleeping,
Dreaming, and Dying:  An Exploration of Consciousness with the Dalai
Lama} (Wisdom Publications:  Boston, 1997).

\item \verb|[11]| F. J. Varela, J.-P. Lachauz, E. Rodrigues and J.
Martinerie, ``The brainweb:  phase synchronization and large-scale
integration,'' \textit{Nature Reviews Neuroscience}, 2 (2001), pp.
229-239.

\item \verb|[12]| Jean Petitot, Francisco J. Varela, Bernard Pachoud
and Jean-Michel Roy, \textit{Naturalizing Phenomenology:  Issues in
Contemporary Phenomenology and Cognitive Science} (Stanford
University Press:  Stanford, 1999)

\item \verb|[13]| Edward L. Robertson, ``An Algebra for Triadic
Relations,'' Technical Report No. 606, Computer Science Department,
Indiana University, Bloomington IN 47404-4101, January 2005; online
at http://www.cs.indiana.edu/ \newline
pub/techreports/TR606.pdf

\item \verb|[14]| E. F. Codd, \textit{The Relational Model for
Database Management:  Version 2} (Addison-Wesley:  Reading MA, 1990)
is the great visionary's most recent and comprehensive statement.

\item \verb|[15]| E. F. Codd, \textit{Cellular Automata} (Academic
Press:  New York, 1968)

\item \verb|[16]| Quinn Tyler Jackson, \textit{Adapting to Babel --
Adaptivity and Context-Sensiti- vity in Parsing:  From
$a^{n}b^{n}c^{n}$ to RNA} (Ibis Publishing:  P.O. Box3083, Plymouth
MA 02361, 2006; for purchasing information, contact Thothic
Technology Partners, LLC, at their website, www.thothic.com).

\end{description}

\end{document}